\pdfoutput=1 %for arXiv submissions
\DeclareSymbolFont{AMSb}{U}{msb}{m}{n}
\documentclass[11pt,noamsfonts,a4paper]{amsart}
\usepackage[left=1.3in, right=1.3in, top=1.3in, bottom=1.3in]{geometry}
\usepackage{mathtools}
\usepackage{braket}
\usepackage{bm} %bold math, like bold mu
\usepackage{enumitem}
\usepackage[charter]{mathdesign}
\usepackage[tracking]{microtype}
\usepackage{tabularx}
\usepackage{etoolbox}
\usepackage{tikz-cd}
\usetikzlibrary{arrows}
\usepackage{verbatim}
\setlist[1]{labelindent=1em}
\setlist[enumerate]{labelsep=*}
\setlist[enumerate,1]{label={\upshape (\roman*)}, ref={\upshape (\roman*)}}

\makeatletter
\newcommand{\leqnomode}{\tagsleft@true\let\veqno\@@leqno}
\newcommand{\reqnomode}{\tagsleft@false\let\veqno\@@eqno}
\makeatother

\usepackage[utf8]{inputenc}

\usepackage[colorlinks=true,urlcolor=blue!70!black,citecolor=blue!60!black,linkcolor=blue!60!black,pdfauthor={Raymond Cheng, Remy van Dobben de Bruyn},pdftitle={Unbounded negativity on rational surfaces in positive characteristic}]{hyperref} %always add last
\tikzset{>={Straight Barb[length=2pt,width=4pt]}, commutative diagrams/arrow style=tikz}

\newcommand{\parref}[1]{{\bf\ref{#1}}}
\newcommand{\parnameref}[2]{{\bf\hyperref[#2]{#1~\ref*{#2}}}}

\makeatletter
\patchcmd{\@settitle}{\uppercasenonmath\@title}{\scshape\Large}{}{}
\patchcmd{\@setauthors}{\MakeUppercase}{\scshape}{}{}
\patchcmd{\maketitle}{\uppercasenonmath\shorttitle}{\shorttitle}{}{}
\makeatother

\setlist{itemsep=0em,topsep=0cm,partopsep=0em,parsep=\lineskip} %if needed

\newcommand{\pointheader}{\refstepcounter{subsection}\textbf{\thesubsection.}}
\newcommand{\point}{\pointheader~}
\newcommand{\tpoint}[1]{\pointheader~{\bf #1. ---}}

\newcommand{\bpoint}[1]{\pointheader~{\bf #1.}}

\newcommand{\psqedsymb}{\(\blacksquare\)}
\renewcommand{\qed}{~\hfill{\psqedsymb}}

\makeatletter
\newcommand*{\coloneqq}{\mathrel{\rlap{%
           \raisebox{0.3ex}{$\m@th\cdot$}}%
           \raisebox{-0.3ex}{$\m@th\cdot$}}%
           =}
\newcommand{\eqqcolon}{=%
           \mathrel{\rlap{%
           \raisebox{0.3ex}{$\m@th\cdot$}}%
           \raisebox{-0.3ex}{$\m@th\cdot$}}}
\makeatother

\newcommand{\bmu}{\pmb\mu}
\DeclareMathOperator{\Char}{char}

\DeclareMathOperator{\pr}{pr}
\DeclareMathOperator{\Spec}{Spec}
\DeclareMathOperator{\NS}{NS}
\newcommand{\punct}[1]{\makebox[0pt][l]{\,#1}} %so that punctuation does not mess with tikz layout

\title[Unbounded negativity on rational surfaces in positive characteristic]
{Unbounded negativity on rational\\ surfaces in positive characteristic}
\author{Raymond Cheng}
\address{Department of Mathematics \\
  Columbia University \\
  2990 Broadway \\
  New York, NY, 10027, United States of America
}
\email{rcheng@math.columbia.edu}
\author{Remy van Dobben de Bruyn}
\address{Department of Mathematics \\ Princeton University \\ Fine Hall \\ Washington Road \\ Princeton, NJ 08544 \\ United States of America}
\address{Institute for Advanced Study \\ 1 Einstein Drive \\ Princeton, NJ 08540 \\ United States of America}
\email{rdobben@math.princeton.edu}
\date{2 March 2021}

\begin{document}
%\leavevmode\vspace{-2em}
%optional negative break for typesetting
\begin{abstract}
  We give explicit blowups of the projective plane in positive characteristic that
  contain smooth rational curves of arbitrarily negative self-intersection, showing
  that the Bounded Negativity Conjecture fails even for rational surfaces in positive
  characteristic.
\end{abstract}
\maketitle
%\vspace{-2em}

\section*{Introduction}
A smooth projective surface \(X\) over an algebraically closed field is said to
have \emph{Bounded Negativity} if there exists a positive integer \(b(X)\) such
that \(C^2 \geq -b(X)\) for any reduced curve \(C \subset X\). A folklore
conjecture, going back to Enriques and discussed in
\cite[Conjecture I.2.1]{Harbourne:Cracow} and \cite[Conjecture 1.1]{BNC}, is the

\textbf{Bounded Negativity Conjecture. ---}
\emph{Any smooth projective surface in characteristic \(0\) has Bounded Negativity.}

The assumption on the characteristic cannot be dropped: if \(C\) is a curve
over \(\bar{\mathbf F}_p\), then the graph \(\Gamma_{F^e} \subseteq C \times C\)
of the \(p^e\)-th power Frobenius endomorphism has self-intersection \(p^e(2-2g)\),
which becomes arbitrarily negative as \(e \to \infty\) when \(g \geq 2\). Nonetheless,
it is conceivable that certain geometric assumptions on the surface may still guarantee
Bounded Negativity in positive characteristic. For instance,
\cite[discussion preceding Example 3.3.3]{RecentDevelopments}
and \cite[Conjecture 2.1.2]{Harbourne-IMPAN} ask whether smooth rational surfaces
over a field of positive characteristic have Bounded Negativity. We give
a negative answer to this question:

\newcounter{intro}
\refstepcounter{intro}\textbf{Main Theorem. ---}\label{theorem}
\emph{Let \(k\) be an algebraically closed field of characteristic \(p > 0\),
let \(m\) be a positive integer invertible in \(k\), and let \(R_m\) be
the blowup of \(\mathbf{P}^2\) along
\[ Z_m \coloneqq \Set{[x_0:x_1:x_2] | x_0^m = x_1^m = x_2^m}. \]
Let \(C_1 = V(x_0+x_1+x_2) \subseteq \mathbf P^2\), and for \(d \geq 1\)
invertible in \(k\), write \(C_d \subseteq \mathbf P^2\) for the image of
\begin{align*}
  \phi_d \colon C_1 &\to \mathbf{P}^2\\
  [x_0:x_1:x_2] &\mapsto [x_0^d:x_1^d:x_2^d].
\end{align*}
If \(dm = p^e-1\) for some positive integer \(e\), then the strict transform
\(\widetilde{C}_d \subseteq R_m\) of \(C_d\) is a smooth rational curve with
\(\widetilde{C}_d^2 = d(3 - m) - 1\). Thus, if \(m > 3\), the rational
surface \(R_m\) does not have Bounded Negativity over \(k\).}

Since \(\mathbf{P}^2\) has Bounded Negativity, this shows that
\cite[Problem 1.2]{BNC:Arrangements} has a negative answer in positive
characteristic:

\textbf{Corollary. ---}
\emph{Bounded Negativity is not a birational property of smooth projective
surfaces in positive characteristic.} \qed

In fact, since every smooth projective surface \(X\) admits a finite morphism
\(X \to \mathbf P^2\), pulling back the blowup \(R_m \to \mathbf P^2\)
gives a blowup \(\widetilde X \to X\) with a finite morphism
\(\widetilde X \to R_m\). Pulling back the curves \(\widetilde C_d\) to
\(\widetilde X\) shows:

\textbf{Corollary. ---}
\emph{If \(X\) is a smooth projective surface over an algebraically closed field
\(k\) of positive characteristic, then there exists a blowup \(\widetilde X \to X\)
such that \(\widetilde X\) does not have Bounded Negativity.} \qed

In \parnameref{\S\!\!}{S:proof}, we give a direct proof of the \hyperref[theorem]{\textbf{Main Theorem}}.
In \parnameref{\S\!\!}{S:line-config}, we realise the plane curves \(C_d\) as norms of line
configuration, thereby deriving equations for them.
In \parnameref{\S\!\!}{S:FerFr}, we view \(R_m\) as an isotrivial family of diagonal
curves over \(C_1\) and relate the curves \(\widetilde C_d\) on \(R_m\) to graphs of
Frobenius morphisms on Fermat curves. Finally, we close in \parnameref{\S\!\!}{S:char-0}
with some questions and remarks towards characteristic zero.

Sections \parref{S:line-config} and \parref{S:FerFr} each give
alternative methods for computing the self-intersections of \(\widetilde C_d\).
Given the simplicity of the formulas for \(\widetilde C_d\) and the many
connections to other well-studied examples, it is surprising that these curves have not
been found before.

\section*{Notation}
Throughout the paper, \(k\) will be an algebraically closed field of
arbitrary characteristic, and \(m\) and \(d\) will denote
positive integers invertible in \(k\). We will use the notation of
the \hyperref[theorem]{\textbf{Main Theorem}} throughout.

\section{Proof of Main Theorem}\label{S:proof}
In this section, fix \(m\) and write \(R \coloneqq R_m\) for the
blowup of \(\mathbf{P}^2\) along \( Z \coloneqq Z_m\).
%optional break for typesetting
The generators \(s_0 = x_1^m-x_2^m\), \(s_1 = x_2^m-x_0^m\),
and \(s_2 = x_0^m-x_1^m\) of the ideal of~\(Z\) give a closed immersion
\(R \hookrightarrow \mathbf P^2 \times \mathbf P^2\).
Since \(s_0+s_1+s_2=0\), one of the \(s_i\) can be eliminated at the expense of
breaking the symmetry in the computations below.

\tpoint{Lemma}\label{about-R}
\emph{The embedding \(R \hookrightarrow \mathbf P^2 \times \mathbf P^2\)
realises \(R\) as the complete intersection
\[
\Set{\big([x_0:x_1:x_2],[y_0:y_1:y_2]\big) \in \mathbf P^2 \times \mathbf P^2
| \begin{array}{c}
y_0+y_1+y_2 = 0\\
x_0^my_0 + x_1^my_1 + x_2^my_2 = 0
\end{array}}
\]
of degrees \((0,1)\) and \((m,1)\) in \(\mathbf{P}^2 \times \mathbf{P}^2\).
In particular, \(K_{R} = \mathcal{O}_{R}(m-3,-1)\).}

\emph{Proof.}
The generators \(s_0, s_1, s_2\) of the ideal of \(Z\) identify
\(R\) as
\[
  \Set{\big([x_0:x_1:x_2], [y_0:y_1:y_2]\big) \in \mathbf{P}^2 \times \mathbf{P}^2
  | \begin{array}{c}
      y_0(x_2^m-x_0^m) = y_1(x_1^m-x_2^m) \\
      y_1(x_0^m-x_1^m) = y_2(x_2^m-x_0^m) \\
      y_2(x_1^m-x_2^m) = y_0(x_0^m-x_1^m)
    \end{array}}.
\]
The relation \(s_0+s_1+s_2 = 0\) shows that \(R\) is contained in the
locus \(y_0+y_1+y_2=0\). The equation \(y_0(x_2^m-x_0^m) = y_1(x_1^m-x_2^m)\)
can be rewritten as \((y_0+y_1)x_2^m = x_0^my_0 + x_1^my_1\), which is
equivalent to \(x_0^my_0+x_1^my_1+x_2^my_2 = 0\) since \(y_0+y_1+y_2=0\).
The same holds for the other two equations by symmetry. The final statement
then follows from the adjunction formula since
\(K_{\mathbf P^2 \times \mathbf P^2} = \mathcal O(-3,-3)\).
\qed

Alternatively, one can observe that the complete intersection of
\parnameref{Lemma}{about-R} maps birationally onto its first factor,
where the fibres are points when \([x_0^m:x_1^m:x_2^m] \neq [1:1:1]\) and
lines otherwise.

\tpoint{Lemma}\label{coord-phi}
\emph{If \(\Char k = p > 0\) and \(dm = p^e-1\) for some positive
integer \(e\), then the map
\(\widetilde\phi_d \colon C_1 \to \mathbf P^2 \times \mathbf P^2\) given by
\[ [x_0:x_1:x_2] \mapsto \big([x_0^d:x_1^d:x_2^d], [x_0:x_1:x_2]\big) \]
lands in \(R\). In particular, it is the unique map lifting
\(\phi_d \colon C_1 \to \mathbf P^2\).
}%

\emph{Proof.}
Since \(x_0+x_1+x_2=0\), the image of \(\widetilde\phi_d\) is contained in the locus
\(y_0+y_1+y_2=0\). Since \(dm = p^e-1\), the expression
\(x_0^my_0+x_1^my_1+x_2^my_2\) pulls back to \(x_0^{p^e}+x_1^{p^e}+x_2^{p^e}\),
which vanishes because the \(p^e\)-th power Frobenius is an endomorphism.
Thus \(\widetilde\phi_d\) is a lift of \(\phi_d\) to \(R\), and it is the unique lift
since the first projection \(\pr_1 \colon R \to \mathbf P^2\) is birational.
\qed

\tpoint{Corollary}\label{self-intersection}
\emph{The map \(\widetilde\phi_d \colon C_1 \to R\) is a closed immersion, whose
image \(\widetilde{C}_d\) is a smooth rational curve in \(R\) with
\(\widetilde{C}_d^2 = d(3-m)-1\).}

\emph{Proof.} The first two statements follow from the coordinate
expression in \parnameref{Lemma}{coord-phi}, since~\(\widetilde\phi_d\) embeds
\(C_1\) linearly into the second factor. The same expression shows that
\(\widetilde\phi_d^*\mathcal{O}_R(a,b) = \mathcal{O}_{C_1}(da+b)\),
so \(K_{R} \cdot \widetilde{C}_d = d(m-3)-1\) by \parnameref{Lemma}{about-R}.
Then the adjunction formula shows that
\reqnomode
\[
  \widetilde{C}_d^2 = -2 - K_{R} \cdot \widetilde{C}_d = d(3-m)-1.\tag*{\psqedsymb}
\]
\leqnomode%
This completes the proof of the \hyperref[theorem]{\textbf{Main Theorem}}. \qed

A consequence of \parnameref{Corollary}{self-intersection} is that the singularities of
\(C_d\) are contained in \(Z\). However, the individual multiplicities are not so
easy to determine. For example, in \parnameref{Lemma}{zeta-Fermat} we will compute the
multiplicity of \(C_d\) at \([1:1:1]\) in terms of point counts on Fermat curves.

\section{Relation with line configurations}\label{S:line-config}
In this section, we observe that the \(d\)-th power maps
\(\pi_d \colon \mathbf{P}^2 \to \mathbf{P}^2\) are finite Galois morphisms such
that \(\pi_d^*C_d\) is the union of the Galois translates of \(C_1\). Thus
the \(C_d\) are norms of line configurations, from which we derive in
\parnameref{Corollary}{cor-norm} a formal product formula for the equation of
the plane curves \(C_d\). In the second half of this section, we observe that
in characteristic \(p > 0\) and for \(q\) a power of \(p\), the curve \(C_{q - 1}\)
comes from a subconfiguration of the set of all \(\mathbf{F}_q\)-rational lines.
This allows us to show in \parnameref{Corollary}{cor-eqn-char-p} that an equation
of \(C_{q-1}\) in this case is the complete homogeneous polynomial of degree \(q-1\).

\bpoint{Power Maps}\label{finite-morphisms}
For any integer \(a \geq 1\) invertible in \(k\), write \(\pi_a\) for the \(a\)-th
power map \(\mathbf P^2 \to \mathbf P^2\). Since \(\pi_a^* Z_m = Z_{am}\),
the map \(\pi_a\) lifts to a finite morphism
\(\widetilde\pi_a \colon R_{am} \to R_m\) given by
\[
\big([x_0:x_1:x_2],[y_0:y_1,y_2]\big) \mapsto \big([x_0^a:x_1^a:x_2^a],[y_0:y_1,y_2]\big).
\]
Since \(a\) is invertible in \(k\), both \(\pi_a\) and
\(\widetilde\pi_a\) are finite Galois with group \(G = \bmu_a^3/\bmu_a\), where
\((\zeta_0,\zeta_1,\zeta_2) \in G\) acts on \(\mathbf P^2\) via
\[
[x_0:x_1:x_2] \mapsto
[\zeta_0x_0:\zeta_1x_1:\zeta_2x_2].
\]
This gives a tower of extensions
\begin{equation*}
\begin{tikzcd}[column sep=.1em,row sep=.8em,arrows=-]
R_4 \ar{rd} & R_6 \ar{d}\ar{rd} & R_9 \ar{d} & \ \scalebox{.6}[1]{$\iddots$} \ar{ld}\\
 & R_2 \ar{rd} & R_3 \ar{d} & \ \scalebox{.6}[1]{$\iddots$}\ar{ld}\\
 & & R_1 & &
\end{tikzcd}
\end{equation*}
indexed by the poset of positive integers invertible in \(k\) under the divisibility relation.

\tpoint{Lemma}\label{totally-split}
\emph{If \(a,d \geq 1\) are invertible in \(k\), then
\begin{enumerate}
\item\label{split-i} The map \(\phi_d \colon C_1 \to \mathbf P^2\) is unramified and
birational onto its image;
\item\label{split-ii} The inverse image \(\pi_a^*C_{ad}\) is totally split into the
\(G\)-translates of \(C_d\).
\end{enumerate}
}

\emph{Proof.}
The Jacobian \((d \cdot x_0^{d-1},d \cdot x_1^{d-1},d \cdot x_2^{d-1})\) of \(\phi_d\) only
vanishes when \(x_0=x_1=x_2=0\), showing that \(\phi_d\) is unramified. Then
the map \(C_1 \to C_d^\nu\) to the normalisation of \(C_d\) is
unramified, hence an isomorphism since it is an \'etale map of smooth
projective rational curves, proving \ref{split-i}.

Since \(\pi_a \circ \phi_d = \phi_{ad}\), part \ref{split-i} shows that \(\pi_a\) maps
\(C_d\) birationally onto its image. This shows that the decomposition group
of \(C_d\) is trivial, so no two \(G\)-translates \(\zeta C_d\) of \(C_d\) coincide
and \(C_{ad}\) is totally split under \(\pi_a\). \qed

\tpoint{Corollary}\label{cor-norm}
\emph{If \(d\) is invertible in \(k\), then the homogeneous ideal of
\(C_d \subseteq \mathbf P^2\) is generated by
\[
f_d \coloneqq N_{\pi_{d,*} \mathcal O_{\mathbf P^2}/\mathcal O_{\mathbf P^2}} \big(x_0^{1/d}+x_1^{1/d}+x_2^{1/d}\big) =
\prod_{\zeta,\zeta' \in \bmu_d} \big(x_0^{1/d}+\zeta x_1^{1/d}+\zeta' x_2^{1/d}\big).
\]
}%
\emph{Proof.}
By \parnameref{Lemma}{totally-split} \ref{split-ii}, the inverse image \(\pi_d^{-1}(C_d)\)
is the union of lines \(\bigcup_{\zeta \in \bmu_d^3/\bmu_d} \zeta C_1\).
The result follows since \(C_1\) is cut out by \(x_0+x_1+x_2=0\).
\qed

\point
In general, the \(f_d\) are complicated symmetric polynomials. However,
in \parnameref{Corollary}{cor-eqn-char-p} we will show that the coefficients of
\(f_{q-1}\) are congruent to \(1\) modulo \(p\) if \(q\) is a power of a prime
\(p\). For example, for \(q = p = 3\), we get
\begin{align*}
  N\big(x_0^{\tfrac{1}{2}}+x_1^{\tfrac{1}{2}}+x_2^{\tfrac{1}{2}}\big) &= \Big(x^{\tfrac{1}{2}}_0+x^{\tfrac{1}{2}}_1+x^{\tfrac{1}{2}}_2\Big)\Big(x_0^{\tfrac{1}{2}}+x_1^{\tfrac{1}{2}}-x_2^{\tfrac{1}{2}}\Big)\Big(x_0^{\tfrac{1}{2}}-x_1^{\tfrac{1}{2}}+x_2^{\tfrac{1}{2}}\Big)\Big(x_0^{\tfrac{1}{2}}-x_1^{\tfrac{1}{2}}-x_2^{\tfrac{1}{2}}\Big) \\
&= x_0^2+x_1^2+x_2^2-2x_0x_1-2x_1x_2-2x_2x_0.\\
&\equiv x_0^2+x_1^2+x_2^2+x_0x_1+x_1x_2+x_2x_0 \pmod 3.
\end{align*}
In the remainder of this section, assume \(\Char k = p > 0\) and let \(q\) be a
power of \(p\).

\bpoint{Finite Field Line Configurations}\label{line-config}
The configuration of \(\mathbf{F}_q\)-rational lines in \(\mathbf{P}^2\) is the
union of the lines \(L_a = \{a_0x_0+a_1x_1+a_2x_2 = 0\}\) indexed by
\(a = [a_0:a_1:a_2] \in \check{\mathbf P}^2(\mathbf F_q)\).
Their union is the divisor in \(\mathbf{P}^2\) cut out by the polynomial
\[
  \det\left(\begin{smallmatrix}
    x_0 & x_1 & x_2 \\ x_0^q & x_1^q & x_2^q \\ x_0^{q^2} & x_1^{q^2} & x_2^{q^2}
  \end{smallmatrix}\right)
  = x_0^q x_1^{q^2} x_2 - x_0^{q^2} x_1^q x_2
  + x_0 x_1^q x_2^{q^2} - x_0 x_1^{q^2} x_2^q
  + x_0^{q^2} x_1 x_2^q - x_0^q x_1 x_2^{q^2},
\]
since the three columns become linearly dependent when \(x_0\), \(x_1\), and \(x_2\)
satisfy a linear relation over \(\mathbf F_q\), and the degree equals
\(q^2+q+1 = |\check{\mathbf P}^2(\mathbf F_q)|\).
%optional break for typesetting
Now \parnameref{Lemma}{totally-split}\ref{split-ii} shows that
\(\pi_{q-1}^*C_{q-1}\) consists of the \(q^2 - 2q + 1\) lines \(L_a\)
with all coordinates of \(a = [a_0:a_1:a_2]\) nonzero.
We can thus derive an equation for \(C_{q-1}\) by extracting factors cutting
out the lines \(L_a\) in which \(a\) has a vanishing coordinate. A neat
description of the result comes from the following polynomial identity,
also observed in \cite[p.\ 90]{RVVZ}:

\tpoint{Lemma}\label{lem-generating-function}
\emph{For any nonnegative integer \(n\), define the polynomials
\[
  g_n \coloneqq \sum_{n_0 + n_1 + n_2 = n} x_0^{n_0} x_1^{n_1} x_2^{n_2}
  \quad\text{and}\quad
  h_n \coloneqq x_0x_1^n - x_0^nx_1 + x_1x_2^n - x_1^nx_2 + x_2x_0^n - x_2^nx_0
\]
in \(\mathbf Z[x_0,x_1,x_2]\). Then \(h_2 = (x_2-x_1)(x_0-x_2)(x_1-x_0)\) and
\(h_n = h_2g_{n-2}\) for \(n \geq 3\).
}

\emph{Proof.}
Let \(G(t) \coloneqq \sum_{n \geq 0} g_n t^n\) and
\(H(t) \coloneqq \sum_{n \geq 0} h_n t^n\) be the generating functions of the
\(g_n\) and \(h_n\), respectively. A standard computation gives
\[
G(t) = \frac{1}{(1 - x_0t)(1 - x_1t)(1 - x_2t)}.
\]
On the other hand, writing
\(h_n = (x_2 - x_1) x_0^n + (x_0 - x_2) x_1^n + (x_1 - x_0) x_2^n\)
gives
\[
  H(t) = \frac{x_2 - x_1}{1 - x_0t} + \frac{x_0 - x_2}{1 - x_1t} + \frac{x_1 - x_0}{1 - x_2t}
  = \frac{(x_2 - x_1)x_2x_1 + (x_0 - x_2)x_0x_2 + (x_1 - x_0)x_0x_1}{(1 - x_0t)(1 - x_1t)(1 - x_2t)} t^2.
\]
The result follows by recognising the numerator as \(h_2\).
\qed

\tpoint{Corollary}\label{cor-eqn-char-p}
\emph{Suppose \(\Char k = p > 0\) and \(q\) is a power of \(p\). Then
\(g_{q-1}\) generates the homogeneous ideal of \(C_{q-1} \subseteq \mathbf P^2\).
In particular, \(f_{q-1} \equiv g_{q-1} \pmod{p}\).
}

\emph{Proof.}
Since \(C_1 = L_{[1:1:1]}\) is among the \(\mathbf{F}_q\)-rational lines of
\parref{line-config} and is not a coordinate axis, the points \([x_0:x_1:x_2]\)
of \(C_1 = V(x_0+x_1+x_2)\) satisfy the equation there divided by \(x_0x_1x_2\):
\[
x_0^{q-1}x_1^{q^2-1}-x_0^{q^2-1}x_1^{q-1}+x_1^{q-1}x_2^{q^2-1}-x_1^{q^2-1}x_2^{q-1}
+x_2^{q-1}x_0^{q^2-1}-x_2^{q^2-1}x_0^{q-1} = 0.
\]
Since \(C_1 \to C_{q-1}\) is \([x_0:x_1:x_2]\mapsto[x_0^{q-1}:x_1^{q-1}:x_2^{q-1}]\),
any \([x_0:x_1:x_2] \in C_{q-1}\) satisfies
\[
  x_0x_1^{q+1}-x_0^{q+1}x_1+x_1x_2^{q+1}-x_1^{q+1}x_2+x_2x_0^{q+1}-x_2^{q+1}x_0 = 0.
\]
So \(h_{q+1}\) vanishes on \(C_{q-1}\), which by
\parnameref{Lemma}{lem-generating-function} equals \((x_2-x_1)(x_0-x_2)(x_1-x_0)g_{q-1}\).
The result follows since \(C_{q-1}\) is not contained in any of the lines
\(\{x_2=x_1\}\), \(\{x_0=x_2\}\), or \(\{x_1=x_0\}\), and \(\deg g_{q-1} = q-1 = \deg C_{q-1}\).
\qed

\bpoint{Negative Curves via Equations}\label{intersections-components}
If \(m > 3\) and \(q\) is a power of \(p\) congruent to \(1\) modulo \(m\), then
the curves \(\widetilde C_d \subseteq R_m\) with \(dm = q^e - 1\) of the
\hyperref[theorem]{\textbf{Main Theorem}} can therefore be obtained by starting
with the very explicit equations
\[
  C_{q^e - 1} = V\left( \sum_{n_0+n_1+n_2=q^e-1} x_0^{n_0} x_1^{n_1} x_2^{n_2} \right) \subseteq \mathbf P^2,
\]
blowing up at \([1:1:1]\), pulling back along
\(\widetilde \pi_m \colon R_m \to R_1\), and taking one of the \(m^2\) isomorphic
components \(\zeta \widetilde C_d\) for \(\zeta \in \bmu_m^3/\bmu_m\). From this
point of view, the self-intersection may be computed as
\[
  \widetilde C_d^2 = \widetilde C_d \cdot \widetilde \pi_m^*(\widetilde C_{q^e-1}) - \sum_{\zeta \neq 1} \widetilde C_d \cdot (\zeta \widetilde C_d) = (2dm-1) - 3(m-1)d = d(3-m)-1,
\]
since the intersection number between \(\widetilde C_d\) and a Galois
translate by \(\zeta = (\zeta_0,\zeta_1,\zeta_2) \in G\setminus \{1\}\) is
\[
\widetilde C_d \cdot (\zeta \widetilde C_d) = \begin{cases}
d, & \zeta_0 = \zeta_1 \text{ or } \zeta_1 = \zeta_2 \text{ or } \zeta_2 = \zeta_0,\\
0, & \text{otherwise}.
\end{cases}
\]
Indeed, \(\zeta\widetilde C_d\) is the image of the morphism
\(\zeta \circ \widetilde \phi_d\) given by
\[
[x_0:x_1:x_2] \mapsto \big([\zeta_0x_0^d:\zeta_1x_1^d:\zeta_2x_2^d],[x_0:x_1:x_2]\big).
\]
Thus, \(\widetilde C_d\) and \(\zeta\widetilde C_d\) only intersect when
\(\zeta\widetilde \phi_d([x_0:x_1:x_2]) = \widetilde \phi_d([x_0:x_1:x_2])\).
At most one of the \(x_i\) can vanish since \(x_0+x_1+x_2=0\), so there are no
intersections when \(\zeta_i \neq \zeta_j\) for \(i \neq j\), and a
single intersection with multiplicity \(d\) at \(V(x_k)\) when \(\zeta_i = \zeta_j\)
and \(\{i,j,k\} = \{0,1,2\}\).

\section{Relation with Fermat varieties and Frobenius morphisms}\label{S:FerFr}

By \parnameref{Lemma}{about-R}, the second projection
\(\pr_2 \colon R_m \to V(y_0+y_1+y_2)\) realises \(R_m\) as the family of diagonal
degree \(m\) curves over \(C_1 \cong \mathbf P^1\) given by
\[
x_0^my_0 + x_1^my_1 + x_2^my_2 = 0.
\]
If \(\Char k = p > 0\) and \(m\) is invertible in \(k\), then the curves
\(\widetilde C_d \subseteq R_m\) for \(dm = p^e-1\) are given by sections
\(\widetilde \phi_d \colon C_1 \to R_m\) of \(\pr_2\).
In this section, we pull back the family \(R_m \to C_1\) and the sections
\(\widetilde \phi_d\) along finite covers of \(C_1\). Pulling back along
covers by Fermat curves allows us to relate the \(\widetilde C_d\) in
\parnameref{Corollary}{Cd-and-graphs} with graphs of Frobenius on products of
Fermat curves. Pulling back along the Frobenius morphism of \(C_1\) allows us
to realise the \(\widetilde C_d\) in \parnameref{Corollary}{rel-Fr} as pullbacks
of a constant section \(\widetilde C_0\) under powers of a horizontal Frobenius
morphism of \(R_m\) over \(C_1\).

\bpoint{Intermediate Surfaces}
For positive integers \(m\) and \(n\) invertible in \(k\) and \(r \in \mathbf N\),
denote by \(R_{m,n,r}\) the normal surface
\[
R_{m,n,r} =
\Set{\big([x_0:x_1:x_2],[y_0:y_1:y_2]\big) \in \mathbf P^2 \times \mathbf P^2
| \begin{array}{c}
y_0^n+y_1^n+y_2^n = 0\\
x_0^my_0^r + x_1^my_1^r + x_2^my_2^r = 0
\end{array}}.
\]
It is smooth if and only if \(m = 1\) or \(r \in\{0,1\}\); in all other cases, the singular
locus \(V(x_0y_0,x_1y_1,x_2y_2)\) consists of the \(3n\) points
\[
\left\{\big([1:0:0],[0:s:t]\big), \big([0:1:0],[s:0:t]\big), \big([0:0:1],[s:t:0]\big)\
\Big|\ s^n+t^n = 0\right\}.
\]
Note that \(R_{m,1,1}\) is none other than the surface \(R_m\) of
\parnameref{Lemma}{about-R}. If \(X_n\) denotes the Fermat curve
\(V(y_0^n+y_1^n+y_2^n) \subseteq \mathbf P^2\) of degree \(n\), then
\(R_{m,n,0}\) coincides with \(X_m \times X_n\). The surfaces \(R_{m,n,r}\)
for \(r > 0\) come with a projection
\[
\pr_2 \colon R_{m,n,r} \to X_n
\]
that is smooth away from the \(3n\) fibres above \(V(y_0y_1y_2) \subseteq X_n\),
and whose singular fibres consist of \(m\) lines meeting at a point.

\bpoint{Generalized Power Maps}\label{pullback}
For positive integers \(a\) and \(b\) invertible in \(k\), define the finite
morphism
\begin{align*}
\pi_{a,b} \colon R_{am,bn,br} &\to R_{m,n,r} \\
\big([x_0:x_1:x_2],[y_0:y_1:y_2]\big) &\mapsto
\big([x_0^a:x_1^a:x_2^a],[y_0^b:y_1^b:y_2^b]\big).
\end{align*}
For \(b=1\) and \(n=r=1\), it coincides with the morphism \(\widetilde \pi_a\)
from \parref{finite-morphisms}. When \(a=1\), these fit into
pullback squares
%\vspace{-.5em}
\begin{equation*}
\begin{tikzcd}
R_{m,bn,br} \ar{r}{\pi_{1,b}}\ar{d}[swap]{\pr_2} & R_{m,n,r}\ar{d}{\pr_2}\\
X_{bn} \ar{r} & X_n\punct{.}
\end{tikzcd}
\end{equation*}
If \(F^e \colon X_n \to X_n\) is the \(p^e\)-th power Frobenius
morphism of \(X_n\), there are pullback squares
%\vspace{-.5em}
\begin{equation*}
\begin{tikzcd}
R_{m,n,p^er} \ar{r}\ar{d}[swap]{\pr_2} & R_{m,n,r} \ar{d}{\pr_2}\\
X_n \ar{r}{F^e} & X_n\punct{,}
\end{tikzcd}
\end{equation*}
so \(R_{m,n,p^er}\) is the Frobenius twist \(R_{m,n,r}^{(e)}\) of \(R_{m,n,r}\)
over \(X_n\). We denote the top map by \(\pi^{(e)}\).

\tpoint{Lemma}\label{birational}
\emph{Let \(m\) and \(n\) be positive integers invertible in \(k\), let \(a,r \in \mathbf Z\),
and assume that \(r\) and \(r+am\) are nonnegative. Then the map
\begin{align*}
  \psi_a \colon \mathbf{P}^2 \times \mathbf{P}^2 & \stackrel{\sim}{\dashrightarrow} \mathbf{P}^2 \times \mathbf{P}^2 \\
  \big([x_0:x_1:x_2],[y_0:y_1:y_2]\big) &\longmapsto
  \big([x_0y_0^a:x_1y_1^a:x_2y_2^a],[y_0:y_1:y_2]\big)
\end{align*}
maps \(R_{m,n,r+am}\) birationally onto \(R_{m,n,r}\).}

\emph{Proof.}
Note that \(\psi_a\) is a birational map with rational inverse \(\psi_{-a}\).
The result follows since \(\psi_a\) takes \(R_{m,n,r+am}\) into \(R_{m,n,r}\)
and \(\psi_{-a}\) does the opposite, and neither surface is contained in the
locus where \(\psi_a\) or \(\psi_{-a}\) is undefined.
\qed

This allows us to relate \(R_{m,m,m}\) and \(X_m \times X_m\):

\tpoint{Corollary}\label{Cd-and-graphs}
\emph{The surfaces \(X_m \times X_m \cong R_{m,m,0}\) and \(R_{m,m,m}\) are birational via
\begin{align*}
\psi \colon X_m \times X_m &\stackrel\sim\dashrightarrow R_{m,m,m}\\
\big([x_0:x_1:x_2],[y_0:y_1:y_2]\big) &\longmapsto
\Big(\big[\tfrac{x_0}{y_0}:\tfrac{x_1}{y_1}:\tfrac{x_2}{y_2}\big],[y_0:y_1:y_2]\Big).
\end{align*}
The composition \(\rho \colon X_m \times X_m \dashrightarrow R_m\) of \(\psi\) with
\(\pi_{1,m}\) is given by
\[
\big([x_0:x_1:x_2],[y_0:y_1:y_2]\big) \longmapsto
\Big(\big[\tfrac{x_0}{y_0}:\tfrac{x_1}{y_1}:\tfrac{x_2}{y_2}\big],[y_0^m:y_1^m:y_2^m]\Big).
\]
If \(\Char k = p > 0\), \(m\) is invertible in \(k\), and \(dm = p^e-1\) for some
positive integer \(e\), then the strict transform of \(\pi_{1,m}^* \widetilde C_d\)
under \(\psi\) is the transpose \(\Gamma_{F^e}^\top\) of the graph of the
\(p^e\)-power Frobenius.
}

\emph{Proof.}
The first statement follows by applying \parnameref{Lemma}{birational} to
\(m=n=r\) and \(a=-1\), and the second is immediate from the definitions. For the
final statement, recall that \(\Gamma_{F^e}^\top\) is given by the section
\(s \colon X_m \to X_m \times X_m\) of \(\pr_2\) given by
\[
[y_0:y_1:y_2] \mapsto \Big(\big[y_0^{p^e}:y_1^{p^e}:y_2^{p^e}\big],[y_0:y_1:y_2]\Big).
\]
By the first pullback square of \parref{pullback}, the curve \(\pi_{1,m}^* \widetilde C_d\)
is the image of the section \(X_m \to R_{m,m,m}\) given by
\[
[y_0:y_1:y_2] \mapsto \big([y_0^{dm}:y_1^{dm}:y_2^{dm}],[y_0:y_1:y_2]\big),
\]
which agrees with \(\psi \circ s\).
\qed\pagebreak

\point
The curves \(\Gamma_{F^e} \subseteq X_m \times X_m\) are the
standard example of curves with unbounded negative self-intersection: the
condition \(m > 3\) of the \hyperref[theorem]{\textbf{Main Theorem}} is exactly
the condition \(g(X_m) > 1\) that makes \(\Gamma_{F^e}^2 = p^e(2-2g)\) negative.
In fact, since \(\Gamma_{F^e}^\top\) passes through \(3m\) of the \(3m^2\) points
of indeterminacy of \(\psi\), resolving the map shows that
\(m^2\,\widetilde C_d^2 = \Gamma_{F^e}^2 - 3m\).

On the other hand, \(R_m \to C_1\) is an isotrivial family of diagonal degree \(m\)
curves that becomes rationally trivialised over the \(m\)-th power cover
\(X_m \to C_1\). Thus, we can also look directly at the pullback
\(\pi^{(e)} \colon R_m^{(e)} \to R_m\) of the Frobenius \(F^e \colon C_1 \to C_1\).
Note that \(R_m^{(e)} = R_{m,1,p^e}\) by \parref{pullback}, so we get:

\tpoint{Corollary}\label{rel-Fr}
\emph{If \(p^e = dm+1\), then \(R_m^{(e)}\) is birational to \(R_m\) via
\begin{align*}
\psi \colon R_m &\stackrel\sim\dashrightarrow R_m^{(e)} \\
\big([x_0:x_1:x_2],[y_0:y_1:y_2]\big) &\longmapsto
\Big(\big[\tfrac{x_0}{y_0^d}:\tfrac{x_1}{y_1^d}:\tfrac{x_2}{y_2^d}\big],[y_0:y_1:y_2]\Big).
\end{align*}
If \(\widetilde \phi_0 \colon C_1 \to R_m\) denotes the constant section
\([y_0:y_1:y_2] \mapsto \big([1:1:1],[y_0:y_1:y_2]\big)\) and
\(\widetilde C_0 \subseteq R_m\) denotes its image, then \(\widetilde C_d\) is the
strict transform of \(\pi^{(e),*} \widetilde C_0\) under \(\psi\).
}%

\emph{Proof.}
The first statement follows from \parnameref{Lemma}{birational} applied to \(n=1\),
\(r = p^e\), and \(a = -d\). For the second, by the second pullback square of \parref{pullback},
the curve \(\pi^{(e),*}\widetilde C_0\) is the image of the constant section \(C_1 \to R_m^{(e)}\)
given by
\[
[y_0:y_1:y_2] \mapsto \big([1:1:1],[y_0:y_1:y_2]\big),
\]
which agrees with \(\psi \circ \widetilde \phi_d\).
\qed

\point
Instead of the transpose \(\Gamma_{F^e}^\top\) of the graph of
\(F^e\colon X_m \to X_m\), one can also look at the negative curves
\(\Gamma_{F^e} \subseteq X_m \times X_m\), which are given by pulling back the
diagonal along the relative Frobenius of \(\pr_2 \colon X_m \times X_m \to X_m\). Their
images under the rational map \(\rho\) of \parnameref{Corollary}{Cd-and-graphs}
are given by the parametrised rational curves
\begin{align*}
C_1 &\to R_m \\
[y_0:y_1:y_2] &\mapsto \Big(\big[y_1^dy_2^d:y_0^dy_2^d:y_0^dy_1^d\big],\big[y_0^{p^e}:y_1^{p^e}:y_2^{p^e}\big]\Big),
\end{align*}
where \(dm = p^e-1\) as usual. These are obtained from the curve \(\widetilde C_0\)
of \parnameref{Corollary}{rel-Fr} by pulling back the strict transform
of \(\widetilde C_0\) under \(\psi^{-1} \colon R_m^{(e)} \dashrightarrow R_m\) along
the relative Frobenius \(F_{R_m/C_1}^e\). Note that the images of these
curves in \(\mathbf P^2\) differ from the curves \(C_d\) by the Cremona transformation
\[
[x_0:x_1:x_2] \mapsto [x_0^{-1}:x_1^{-1}:x_2^{-1}].
\]
Finally, we relate the multiplicity of \(C_d\) at \([1:1:1]\) to point counts
on the Fermat curve \(X_m\) if \(dm = p^e-1\) for some positive integer \(e\).

\tpoint{Lemma}\label{zeta-Fermat}
\emph{If \(dm = p^e-1\), then a point \(x \in C_1\) maps to \([1:1:1]\)
in \(C_d\) if and only if there exists \(y \in X_m(\mathbf F_{p^e})\)
with nonzero coordinates mapping to \(x\) under the \(m\)-th power map
\(X_m \to C_1\). In particular,
\[
\operatorname{mult}_{[1:1:1]} C_d = \frac{|X_m(\mathbf F_{p^e})|-3m}{m^2}.
\]
}%
\emph{Proof.}
The first statement follows since \(X_m \to C_1\) is surjective and a point
\(y = [y_0:y_1:y_2]\) on \(X_m\) with nonzero coordinates maps to \([1:1:1]\)
under the \((p^e-1)\)-st power map \(X_m \to C_d\) if and only if
\(y \in X_m(\mathbf F_{p^e})\).
%optional break for typesetting

For the second statement, note that
\(\operatorname{mult}_{[1:1:1]} C_d\) equals the number of preimages of \([1:1:1]\)
in \(C_1\), since \(\widetilde \phi_d \colon C_1 \to \widetilde C_d\) is an
isomorphism by \parnameref{Corollary}{self-intersection}. The result now
follows since \(X_m \to C_1\) is finite \'etale of degree \(m^2\) away from the
coordinate axes, so each point \(x \in C_1 \setminus V(x_0x_1x_2)\) has exactly
\(m^2\) preimages in \(X_m\).
\qed

For example, if \(p^\nu \equiv -1 \pmod m\) for some positive integer \(\nu\),
then
\[
|X_m(\mathbf F_{p^e})| = 1 - \frac{(m-1)(m-2)}{2} p^{e/2} + p^e
\]
whenever \(p^e \equiv 1 \pmod m\) \cite[Lem.\ 3.3]{ShiodaKatsura}.

\section{Remarks towards characteristic \texorpdfstring{\(0\)}{0}}\label{S:char-0}

\point
Although the Bounded Negativity Conjecture is currently still open in
characteristic \(0\), the Weak Bounded Negativity Conjecture is known \cite{Hao}:
for any smooth projective complex surface \(X\) and any \(g \in \mathbf N\),
there exists a constant \(b(X,g)\) such that \(C^2 \geq -b(X,g)\) for every
reduced curve \(C = \sum_i C_i\) whose components \(C_i\) have geometric
genus at most \(g\).

Our examples in the \hyperref[theorem]{\textbf{Main Theorem}} certainly violate
this, and, as we now verify, arise from the failure of the logarithmic
Bogomolov--Miyaoka--Yau inequality for the pair \((R_m,\widetilde C_d)\) when
\(d\) is large with respect to \(m\). In the next three paragraphs, assume
\(\Char k = p > 0\) and \(dm = p^e - 1\). To ease notation, write
\((R,\widetilde C)\) for \((R_m, \widetilde C_d)\).
We will use logarithmic sheaves of differentials; see for example
\cite[\S2]{EsnaultViehweg}.

\tpoint{Lemma}\label{chern-numbers}
\emph{The Chern numbers of the pair \((R, \widetilde C)\) are
  \begin{align*}
    c_1^2(R,\widetilde C)
    & \coloneqq c_1^2\big(\Omega^1_{R}(\log \widetilde C)\big)
      = d(m-3) - m^2 + 6, \\
    c_2(R,\widetilde C)
    & \coloneqq c_2\big(\Omega^1_{R}(\log \widetilde C)\big)
      = m^2 + 1.
  \end{align*}
  In particular, the Chern slopes \(c_1^2(R,\widetilde C)/c_2(R,\widetilde C)\)
  are unbounded for fixed \(m\) and growing \(d\).}

\emph{Proof.}
The logarithmic sheaf of differentials fit into a short exact sequence
\[
  0 \to
  \Omega^1_{R} \to
  \Omega^1_{R}(\log \widetilde C) \to
  \mathcal{O}_{\widetilde C} \to
  0,
\]
so \(c_1^2(R,\widetilde C) = (K_{R} + \widetilde C)^2\) and
\(c_2(R,\widetilde C) = c_2(\Omega^1_{R}) + \widetilde C(K_{R} + \widetilde C)\).
Since \(R\) is the blowup of \(\mathbf P^2\) in \(m^2\) points, we get
\(K_R^2=9-m^2\) and \(c_2(\Omega^1_R) = 3+m^2\), so the result follows from the
computations of the intersection numbers in \parnameref{Corollary}{self-intersection}.
\qed

\tpoint{Lemma}\label{KC-pseff}
\emph{If \(m > 3\) and \(d\) is such that
\[
  \chi\big(2(K_{R} + \widetilde C)\big) = d(m-3) - m^2 + 5 > 0
\]
then \(\mathrm{H}^0(R, 2(K_{R} + \widetilde C)) \neq 0\).
In particular, \(K_{R} + \widetilde C\) is pseudoeffective.}

\emph{Proof.}
The Euler characteristic statement follows from Riemann--Roch, so it remains
to show that \(\mathrm{H}^0\big(R,2(K_{R} + \widetilde C)\big) \neq 0\) once
\(\chi\big(2(K_{R} + \widetilde C)\big) > 0\).
But \(\mathrm{H}^2\big(R,2(K_R+\widetilde C)\big) = \mathrm{H}^0(R,-K_R-2\widetilde C)^\vee\),
and the latter vanishes since \(\widetilde C\) is effective and
\(-K_R = \mathcal{O}_R(3-m,1)\) by \parnameref{Lemma}{about-R}.
\qed

For \(d\) large with respect to \(m\), this shows that
\((R, \widetilde C)\) falls into the final case considered in \cite[\S 1.2, Case 2]{Hao},
and that the failure of Weak Bounded Negativity stems from the failure of the
logarithmic Bogomolov--Miyaoka--Yau inequality:

\tpoint{Corollary}
\emph{If \(m > 3\) and
\(d > \frac{5m^2 - 2}{m - 3}\), then \(K_{R} + \widetilde C\) is pseudoeffective
and
\[ c_1^2(R,\widetilde C)/c_2(R,\widetilde C) > 4. \]
Moreover, the pair \((R, \widetilde C)\) does not lift to the second Witt
vectors \(W_2(k)\).}

\emph{Proof.}
The first part follows from \parnameref{Lemma}{chern-numbers} and
\parnameref{Lemma}{KC-pseff}. The final statement follows from
\cite[Proposition 4.3]{Langer}, since \((R, \widetilde C)\) violates the
logarithmic Bogomolov--Miyaoka--Yau inequality.
\qed

\point
On the other hand, the surface \(R_m\) itself does lift to characteristic \(0\).
This gives new examples of surfaces \(X \to \Spec \mathbf Z\) such that almost all
special fibres \(X_{\bar{\mathbf F}_p}\) (namely those with \(p \nmid m\)) violate
bounded negativity. The same property holds for the square \(C \times C\) of a curve
\(C \to \Spec \mathbf Z\) of genus \(\geq 2\), which is the classical counterexample
to bounded negativity in positive characteristic. However, the rational surface
\(X = R_m\) has the additional property that the specialisation maps
\(\NS(X_{\bar{\mathbf Q}}) \to \NS(X_{\bar{\mathbf F}_p})\) are isomorphisms for
every prime \(p \nmid m\).

\tpoint{Question}
\emph{Is it possible to determine the effective cone of \(R_m\) for some
\(m \geq 4\)? How does it depend on the characteristic of \(k\)?
}

For example, the curves in \(\mathbf P^2\) cut out by the polynomials
\(g_{m-1}\) of \parnameref{Lemma}{lem-generating-function} are smooth of genus
\(\tfrac{(m-2)(m-3)}{2}\) in characteristic \(0\) \cite[Thm.\ 1]{RVVZ}, and the
equation \(g_{m-1}h_2 = h_{m+1}\) shows that
\(V(g_{m-1}) \cup V(x_0-x_1) \cup V(x_1-x_2) \cup V(x_2-x_0)\) contains
\[
Z_m' \coloneqq V\left(\begin{array}{c}
x_0x_1^{m+1}-x_0^{m+1}x_1,\\
x_1x_2^{m+1}-x_1^{m+1}x_2,\\
x_2x_0^{m+1}-x_2^{m+1}x_0\end{array}\right)
= Z_m \cup \left\{[s:t:0],[s:0:t],[0:s:t]\ \big|\ s^m=t^m\right\}.
\]
Since \(V(g_{m-1})\) has self-intersection \((m-1)^2\) and passes through the
\(m^2-3m+2\) points of \(Z_m\) whose coordinates are pairwise distinct,
its strict transform on \(R_m\) has self-intersection \(m-1\).
On the further blowup \(R'_m\) of \(\mathbf P^2\) in \(Z'_m\), the strict
transform has self-intersection \(-2m+2\), but unlike the situation described
in \parref{intersections-components}, there does not appear to be an obvious way to
produce infinitely many negative curves on a \emph{single} rational surface this way.

When \(m = p^e\) for some prime \(p\), the specialisation to characteristic \(p\)
collapses \(Z_m\) onto the point \([1:1:1]\), and the smooth curve \(V(g_{m-1})\)
becomes a rational curve that is highly singular at \([1:1:1]\). Even though these
curves are not negative yet (see \parref{intersections-components}), taking different
values of \(e\) does give infinitely many curves on the same rational surface.

\point
As far as we are aware, all known counterexamples to bounded negativity on
a smooth projective surface \(X\) over an algebraically closed field \(k\)
of characteristic \(p > 0\) consist of a family \(C_i\) of curves on \(X\)
for which there exist constants \(a,b\) such that \(C_i^2 = ap^i + b\) for all
\(i \in \mathbf N\).

\tpoint{Question}
\emph{If \(X\) is a surface over an algebraically closed field \(k\) of
characteristic \(p > 0\), is there a finite set
\(\{(a_i,b_i) \in \mathbf Q^2\ |\ i \in I\}\) such that all integral curves
\(C \subseteq X\) with \(C^2 < 0\) satisfy
\[
C^2 = a_ip^e+b_i
\]
for some positive integer \(e\) and some \(i \in I\)? If not, is there some
other way in which the self-intersections of negative curves on \(X\) are
``not too scattered''?
}

We can also consider the following uniform version:

\tpoint{Question}\label{q-uniform}
\emph{If \(X \to S\) is a smooth projective surface over a finitely generated
integral base scheme \(S\), does there exist
a finite set \(\{(a_i,b_i) \in \mathbf Q^2\ |\ i \in I\}\) such that every
geometrically integral curve \(C \subseteq X_s\) of negative self-intersection
in a fibre \(X_s\) with \(\Char \kappa(s) > 0\) satisfies
\[
C^2 = a_ip^e+b_i
\]
for some positive integer \(e\) and some \(i \in I\), where \(p = \Char \kappa(s)\)?
}

For example, for the surfaces \(R_m \to \Spec\mathbf Z [1/m]\) and the curves
\(\widetilde C_d\) of the \hyperref[theorem]{\textbf{Main Theorem}}, we may take
\(a=\tfrac{3-m}{m}\) and \(b=\tfrac{-3}{m}\), which do not depend on
the characteristic of \(\kappa(s)\).

\point
Despite the failure of bounded negativity in positive characteristic, a
positive answer to \parnameref{Question}{q-uniform} still implies bounded
negativity in characteristic \(0\) via reduction modulo primes.
Indeed, the minimum \(b_{\text{min}} = \min\{b_i\ |\ i \in I\}\) is a lower
bound for the self-intersection \(C^2\) of a geometrically integral curve \(C\)
in the generic fibre, since the specialisations \(C_s\) of \(C\) satisfy
\(C_s^2=C^2\) for all \(s \in S\) and remain geometrically integral for \(s\)
in a dense open set \(U \subseteq S\), and
\[
\bigcap_{p \in P} \big\{a_ip^e+b_i\ \big|\ i \in I, e \in \mathbf Z_{>0}\big\}
\subseteq [b_{\text{min}},\infty)
\]
for any infinite set of primes \(P\). Thus, \parnameref{Question}{q-uniform}
is a natural analogue of the Bounded Negativity Conjecture in positive
characteristic.

\section*{Acknowledgements}
{\small
We thank Johan de Jong, Joaqu\'in Moraga, Takumi Murayama,
Will Sawin, and John Sheridan for helpful discussions.
RvDdB was partly supported by the Oswald Veblen Fund at the
Institute for Advanced Study.
}

\bibliographystyle{alphaurledit}
\bibliography{main}
\end{document}